\documentclass[11pt]{article}

\usepackage{amssymb,amsmath,amstext,amscd,amsthm,latexsym}

\usepackage[dvipdfmx, bookmarkstype=toc, colorlinks=true, linkcolor=cyan, pdfborder={0 0 0}, bookmarks=true, bookmarksnumbered=true]{hyperref}

\DeclareMathOperator*{\regprod}{\mathchoice%
{\ooalign{\hbox{$\displaystyle\prod$}\crcr\hbox{$\displaystyle\coprod$}}}
{\ooalign{\hbox{$\textstyle\prod$}\crcr\hbox{$\textstyle\coprod$}}}
{\ooalign{\hbox{$\scriptstyle\prod$}\crcr\hbox{$\scriptstyle\coprod$}}}
{\ooalign{\hbox{$\scriptscriptstyle\prod$}\crcr\hbox{$\scriptscriptstyle\coprod$}}}
}

\numberwithin{equation}{section}

\theoremstyle{plain}

\newtheorem{theorem}{Theorem}

\newtheorem{remark}{Remark}

\begin{document}

\title{Determinants of Riemann operators on Quillen's higher $K$--groups: periodicity}

\author{Nobushige Kurokawa\footnote{Department of Mathematics, Tokyo Institute of Technology} \and Hidekazu Tanaka\footnote{6-15-11-202 Otsuka, Bunkyo-ku, Tokyo}}

\date{September 28, 2022}



\maketitle

\begin{abstract}
In a previous paper \cite{KT} we introduced determinant of the Riemann operator on Quillen's higher $K$--groups of the integer ring of an algebraic number field $K$. We showed that the determinant expresses essentially the inverse of the so called gamma factor of Dedekind zeta function of $K$. 

Here we study the periodicity of determinant. This comes from the famous ``periodicity'' of higher $K$ groups. This periodicity is analogous to Euler's periodicity of gamma function $\Gamma(x+1)=x\Gamma(x)$. We investigate the ``reflection formula'' corresponding to Euler's reflection formula $\Gamma(x)\Gamma(1-x)=\frac{\pi}{\sin(\pi x)}$ also.
\end{abstract}

\section*{Introduction} Let $A$ be the integer ring of an algebraic number field $K$ that is a finite extension field of the rational number field ${\mathbb Q}$. We denote by $K_{n}(A)$ the $K$--group constructed by Quillen, and we introduce the Riemann operator $\mathcal{R}$ on $K_{n}(A)$ by
\[
\mathcal{R} | K_{n}(A) = \frac{1-n}{2}.
\]
In this paper we study the determinant
\begin{align*}
G_{K}(s) &= {\rm det} ((s I-\mathcal{R}) | \bigoplus_{n>1} K_{n}(A)_{\mathbb{C}} )\\
&= \regprod_{n>1} (s-\frac{1-n}{2})^{{\rm rank} (K_{n}(A))}\\
&= \regprod_{n>1} (\frac{n-1}{2}+s)^{{\rm rank} (K_{n}(A))}
\end{align*}
as the regularized product, which means that
\[
G_{K}(s) = \exp( -\frac{\partial}{\partial w} \varphi(w,s)\biggl|_{w=0} )
\]
for
\[
\varphi(w,s)=\sum_{n>1} {\rm rank}(K_{n}(A)) (\frac{n-1}{2}+s)^{-w}.
\]
The calculation in \cite{KT} using the result of Borel \cite{B} given by for $n > 1$ as
\[
{\rm rank}(K_{n}(A)) = \left\{
\begin{array}{ccc}
r_{1}+r_{2} & {\rm if} & n \equiv 1 \; {\mod 4},\\
r_{2} & {\rm if} & n \equiv 3 \; {\mod 4},\\
0 & {\rm if} & {\rm otherwise}
\end{array}
\right.
\]
implies that
\begin{align*}
G_{K}(s) &= \biggl( 2^{-\frac{s}{2}} \frac{\sqrt{\pi}}{\Gamma(\frac{s}{2}+1)} \biggl)^{r_{1}+r_{2}} \biggl( 2^{-\frac{s}{2}} \frac{\sqrt{2\pi}}{\Gamma(\frac{s+1}{2})} \biggl)^{r_{2}}\\
&= s^{-r_{1}-r_{2}} \Gamma_{{\mathbb R}}(s)^{-r_{1}} \Gamma_{{\mathbb C}}(s)^{-r_{2}} (2\pi)^{-\frac{[K:{\mathbb Q}]}{2}s} C(K),
\end{align*}
where $r_{1}$ (resp. $r_{2}$) is the number of real (resp. complex) places of $K$ and $[K:{\mathbb Q}]=r_{1}+2r_{2}$ with
\[
C(K) = (2 \sqrt{\pi})^{r_{1}} (2\sqrt{2\pi})^{r_{2}}.
\]
The usual notation $\Gamma_{\mathbb{R}}(s)$ and $\Gamma_{\mathbb{C}}(s)$
are defined as
\[
\Gamma_{\mathbb{R}}(s) = \Gamma(\frac{s}{2})\pi^{-\frac{s}{2}}
\]
and
\[
\Gamma_{\mathbb{C}}(s) = \Gamma(s)2(2\pi)^{-s} = \Gamma_{\mathbb{R}}(s) \Gamma_{\mathbb{R}}(s+1).
\]
We review the calculation in \S 1 below. We prove the periodicity of $G_{K}(s)$ as
\begin{theorem}[periodicity] 
\[
G_{K}(s) = G_{K}(s+2) (s+2)^{r_{1}+r_{2}} (s+1)^{r_{2}}.
\]
\end{theorem}
We present two proofs. The first proof uses the explicit calculation of $G_{K}(s)$ given in \cite{KT}. The second proof is coming from the periodicity in $K$--theory
\[
{\rm rank}(K_{n}(A)) = {\rm rank}(K_{n+4}(A)).
\]
This periodicity is an analogue of Euler's periodicity $\Gamma(s+1)=\Gamma(s)s$.

Next, we study the analogue of Euler's reflection formula
\[
\Gamma(s)\Gamma(1-s)=\frac{\pi}{\sin(\pi s)},
\]
and prove the following
\begin{theorem}
(1)[reflection formula]
\[
G_{K}(s)G_{K}(-s) = \biggl( \frac{2}{s} \sin(\frac{\pi s}{2}) \biggl)^{r_{1}} \biggl( \frac{2}{s} \sin(\pi s) \biggl)^{r_{2}}.
\] 
(2) Especially, $G_{K}(s)G_{K}(-s) \in \overline{{\mathbb Q}}$ for $s \in {\mathbb Q}^{\times}$.
\end{theorem}

\section{Calculations of regularized products} We recall the explicit calculations of $G_{K}(s)$ following \cite{KT}. The result of Borel \cite{B} gives
\begin{align*}
\varphi(w,s)&=\sum_{n>1} {\rm rank}(K_{n}(A)) (\frac{n-1}{2}+s)^{-w}\\
&=(r_{1}+r_{2}) \varphi_{1}(w,s) + r_{2} \varphi_{2}(w,s),
\end{align*}
where
\[
\varphi_{1}(w,s)=\sum_{n > 1 \atop n \equiv 1 \; {\rm mod} \; 4} (\frac{n-1}{2}+s)^{-w}
\]
and
\[
\varphi_{2}(w,s)=\sum_{n \equiv 3 \; {\rm mod} \; 4} (\frac{n-1}{2}+s)^{-w}.
\]
Then Lerch's formula implies
\begin{align*}
\regprod_{n>1 \atop n \equiv 1 \; {\rm mod} \; 4} (\frac{n-1}{2}+s) &= \exp(-\frac{\partial}{\partial w} \varphi_{1}(w,s)\biggl|_{w=0})\\
&= 2^{-\frac{s}{2}} \frac{\sqrt{\pi}}{\Gamma(\frac{s}{2}+1)}\\
&= s^{-1} \Gamma_{\mathbb{R}}(s)^{-1} (\sqrt{2\pi})^{-s} 2 \sqrt{\pi}
\end{align*}
and
\begin{align*}
\regprod_{n \equiv 3 \; {\rm mod} \; 4} (\frac{n-1}{2}+s) &= \exp(-\frac{\partial}{\partial w} \varphi_{2}(w,s)\biggl|_{w=0})\\
&= 2^{-\frac{s}{2}} \frac{\sqrt{2\pi}}{\Gamma(\frac{s+1}{2})}\\
&= \Gamma_{\mathbb{R}}(s+1)^{-1} (\sqrt{2\pi})^{-s} \sqrt{2}.
\end{align*}
Hence we obtain
\[
G_{K}(s) = \biggl( 2^{-\frac{s}{2}} \frac{\sqrt{\pi}}{\Gamma(\frac{s}{2}+1)} \biggl)^{r_{1}+r_{2}} \biggl( 2^{-\frac{s}{2}} \frac{\sqrt{2\pi}}{\Gamma(\frac{s+1}{2})} \biggl)^{r_{2}}
\]
as in \cite{KT}. We notice that in \cite{KT} we calculated
\[
{\rm det} ((s I-\mathcal{R}) | \bigoplus_{n\geq0} K_{n}(A)_{\mathbb{C}})
=(s-\frac{1}{2}) s^{r_{1}+r_{2}-1} {\rm det} ((s I-\mathcal{R}) | \bigoplus_{n>1} K_{n}(A)_{\mathbb{C}}).
\]
We refer to \cite{KT} for detailed calculation and we refer to Deninger \cite{D} and Manin \cite{M} concerning regularized products in general.

\section{Periodicity: Proofs of Theorem 1}

\begin{proof}[First Proof] From the explicit formula
\[
G_{K}(s) = \biggl( 2^{-\frac{s}{2}} \frac{\sqrt{\pi}}{\Gamma(\frac{s}{2}+1)} \biggl)^{r_{1}+r_{2}} \biggl( 2^{-\frac{s}{2}} \frac{\sqrt{2\pi}}{\Gamma(\frac{s+1}{2})} \biggl)^{r_{2}}
\]
we get
\begin{align*}
G_{K}(s+2) &= \biggl( 2^{-\frac{s+2}{2}} \frac{\sqrt{\pi}}{\Gamma(\frac{s+2}{2}+1)} \biggl)^{r_{1}+r_{2}} \biggl( 2^{-\frac{s+2}{2}} \frac{\sqrt{2\pi}}{\Gamma(\frac{s+3}{2})} \biggl)^{r_{2}}\\
&= \biggl( 2^{-\frac{s}{2}} \frac{\sqrt{\pi}}{2(\frac{s+2}{2})\Gamma(\frac{s+2}{2})} \biggl)^{r_{1}+r_{2}} \biggl( 2^{-\frac{s}{2}} \frac{\sqrt{2\pi}}{2(\frac{s+1}{2})\Gamma(\frac{s+1}{2})} \biggl)^{r_{2}}\\
&= (s+2)^{-(r_{1}+r_{2})} (s+1)^{-r_{2}} G_{K}(s).
\end{align*}

\end{proof}

\begin{proof}[Second Proof] From
\[
G_{K}(s) = \regprod_{n>1} (\frac{n-1}{2}+s)^{{\rm rank}(K_{n}(A))}
\]
we get
\begin{align*}
G_{K}(s+2) &= \regprod_{n>1} (\frac{n-1}{2}+s+2)^{{\rm rank}(K_{n}(A))}\\
&= \regprod_{n>1} (\frac{(n+4)-1}{2}+s)^{{\rm rank}(K_{n}(A))}.
\end{align*}
Hence using the periodicity
\[
{\rm rank} (K_{n}(A)) = {\rm rank} (K_{n+4}(A))
\]
we obtain
\begin{align*}
G_{K}(s+2) &= \regprod_{n>1} (\frac{(n+4)-1}{2}+s)^{{\rm rank}(K_{n+4}(A))}\\
&= \regprod_{n>5} (\frac{n-1}{2}+s)^{{\rm rank}(K_{n}(A))}.
\end{align*}
Thus we see that
\begin{align*}
G_{K}(s) &= G_{K}(s+2) (\frac{3-1}{2}+s)^{{\rm rank}(K_{3}(A))} (\frac{5-1}{2}+s)^{{\rm rank}(K_{5}(A))} \\
&= G_{K}(s+2) (s+1)^{r_{2}} (s+2)^{r_{1}+r_{2}},
\end{align*}
since ${\rm rank}(K_{3}(A))=r_{2}$ and ${\rm rank}(K_{5}(A))=r_{1}+r_{2}$.
\end{proof}

\section{Reflection formula: Proof of Theorem 2}

\begin{proof}[Proof of Theorem 2] We prove Theorem 2 by using the explicit formula
\[
G_{K}(s)=\biggl( 2^{-\frac{s}{2}} \frac{\sqrt{\pi}}{\Gamma(\frac{s}{2}+1)} \biggl)^{r_{1}+r_{2}} \biggl( 2^{-\frac{s}{2}} \frac{\sqrt{2\pi}}{\Gamma(\frac{s+1}{2})} \biggl)^{r_{2}}.
\]
We calculate 
\begin{align*}
G_{K}(s)G_{K}(-s) &= \biggl( 2^{-\frac{s}{2}} \frac{\sqrt{\pi}}{\Gamma(\frac{s}{2}+1)} 2^{\frac{s}{2}} \frac{\sqrt{\pi}}{\Gamma(-\frac{s}{2}+1)} \biggl)^{r_{1}+r_{2}} \biggl( 2^{-\frac{s}{2}} \frac{\sqrt{2\pi}}{\Gamma(\frac{s+1}{2})} 2^{\frac{s}{2}} \frac{\sqrt{2\pi}}{\Gamma(\frac{-s+1}{2})} \biggl)^{r_{2}}\\
&= \biggl( \frac{\pi}{\Gamma(\frac{s}{2})\Gamma(1-\frac{s}{2})\frac{s}{2}} \biggl)^{r_{1}+r_{2}} \biggl( \frac{2\pi}{\Gamma(\frac{s+1}{2})\Gamma(1-\frac{s+1}{2})} \biggl)^{r_{2}}\\
&= \biggl( \frac{\sin(\frac{\pi s}{2})}{\frac{s}{2}} \bigg)^{r_{1}+r_{2}} \biggl( 2 \sin(\frac{\pi(s+1)}{2}) \biggl)^{r_{2}}\\
&= \biggl( \frac{2}{s} \sin(\frac{\pi s}{2}) \bigg)^{r_{1}+r_{2}} \biggl( 2 \cos(\frac{\pi s}{2}) \biggl)^{r_{2}}\\
&= \biggl( \frac{2}{s} \sin(\frac{\pi s}{2}) \bigg)^{r_{1}} \biggl( \frac{2}{s} \sin(\frac{\pi s}{2}) 2 \cos(\frac{\pi s}{2}) \biggl)^{r_{2}}\\
&= \biggl( \frac{2}{s} \sin(\frac{\pi s}{2}) \bigg)^{r_{1}} \biggl( \frac{2}{s} \sin(\pi s) \bigg)^{r_{2}},
\end{align*}
which gives (1). In particular, (2) follows from (1), since $\frac{\sin(\pi a)}{a} \in \overline{{\mathbb Q}}$ for $a \in {\mathbb Q}^{\times}$.

\begin{remark}
\[
G_{K}(0)=\pi^{\frac{r_{1}+r_{2}}{2}} 2^{\frac{r_{2}}{2}} \notin \overline{{\mathbb Q}}.
\]
\end{remark}

\end{proof}

\end{document}